\documentclass[12pt]{article}       

\begin{document}
\title{Non-commutative low dimension spaces and superspaces associated with contracted quantum groups and supergroups}
\author{N.A.Gromov, D.B.Efimov, I.V.Kostyakov, V.V.Kuratov  \\
Department of Mathematics, \\ Institute of Mathematics and Mechanics, \\
Ural Division, Russian Academy of Sciences, \\
Chernova st., 3a, Syktyvkar, 167982, Russia \\
e-mail: gromov@dm.komisc.ru}

\maketitle

\begin{abstract}
Quantum planes which correspond to all one parameter solutions of QYBE for the two-dimensional case of GL-groups are summarized and their geometrical interpretations are given. It is shown that the quantum dual plane is associated with an exotic solution of QYBE and the well-known quantum $h$-plane may be regarded as the quantum analog of the flag (or fiber) plane.
Contractions of the quantum supergroup $ GL_q(1|2)$ and corresponding quantum superspace $ C_q(1|2)$ are considered in Cartesian basis. The contracted quantum superspace $ C_h(1|2;\iota)$
is interpreted as the non-commutative analog of the superspace with the fiber odd part.
\end{abstract}

\section{Introduction}

In the quantum group theory solutions of Yang-Baxter equation (YBE) have
an inportant meaning. In particular, the Poincar\'e-Birkhoff-Witt
theorem is hold for such quantum group, which mean that the dimension of the space of homogeneous polynomials is the same as in the commutative case. There are several known solutions of YBE. Some of them have clear
group theoretical and geometrical meaning of non-commutative (quantum) analogs of simple Lie groups, algebras and corresponding spaces \cite{FRT-89}, but the other have not. Moreover a contracted quantum group
may have such commutation relations for generators which can not be written
in the form of RTT-relations (\ref{2}) with some matrix $R$ in spite of
the fact that the initial quantum group is defined by the solution $R_q$
of YBE. Therefore it is important to  analyse the quantum group theory
from geometrical point of view and try to find the commutative analogs
especially in the case of contraction. In this paper we perform such
analysis for the simplest case of the low dimension quantum (super)spaces
and corresponding quantum (super)groups associated with all one parameter
$4 \times 4$  solutions of YBE and for the case of quantum superspace
$C_q(1|2).$

The paper is organized as follows. The basic concepts of the quantum group theory are brifly remined in Section 2. Quantum planes and quantum groups
generated by the standard and exotic solutions of QYBE as well as their
contractions are regarded in Section 3. Section 4 is devoted to the standard and exotic solutions of GYBE and contractions of the superplanes.
Quantum superspace $C_q(1|2)$ and its contractions are regarded in Section 5. Our results are summarized in Conclusion.

\section{Solutions of YBE, quantum groups and  spaces}      
            
  According to the well-known theory \cite{FRT-89} with each
  solution $R_q$ of quantum Yang-Baxter equation
  \begin{equation} 
  R_{12}R_{13}R_{23}=R_{23}R_{13}R_{12}
  \label{1}  
  \end{equation} 
is connected a function algebra on quantum group $G_{q}$
(or simply quantum group), which generators $T=(t_{ij}), i,j=1,\ldots,N$ are subject of  the commutation relations
  \begin{equation} 
R_qT_1T_2=T_2T_1R_q,  
  \label{2}  
  \end{equation} 
and a function algebra on quantum vector space $C^N_R$  
(or simply quantum vector space), which generators 
$X=(x_i)$ commute (in the case of GL-groups) as follows 
  \begin{equation} 
(\hat{R}_q - qI)X \otimes X =0,  
  \label{3}  
  \end{equation} 
where $T_1=T \otimes I, T_2=I \otimes T, \hat{R}_q=PR_q,
Pu\otimes v=v\otimes u $ for any $u,v\in C^N. $ 
A co-action of $G_{q}$ on $C^N_R$ is given by
  \begin{equation} 
\delta(X)=T\stackrel{.}\otimes X, \quad
\delta(x_i)=\sum_{k=1}^Nt_{ik}\otimes x_k.  
  \label{4}  
  \end{equation} 

For $N=2$ one-parameter solutions of YBE are obtained 
\cite{H-92} in the form
\begin{equation}
R_q=\left(\begin{array}{cccc}
q &0& 0 &0 \\
0& 1 &0  &0 \\
0& \lambda & 1 &0 \\
0&0 & 0& \zeta 
\end{array} \right),  
   \label{5}  
  \end{equation} 
where
$\lambda =q-q^{-1}, $
$ \zeta=q $ correspond to the standard solution of QYBE for 
quantum group $GL_q(2);$
$\zeta =-q^{-1} $ --- to the exotic solution of QYBE for $ \tilde{GL}_q(2); $
$ \zeta =q^{-1} $ --- to the standard solution of GYBE  for quantum supergroup $ GL_q(1|1) $ and
$ \zeta =-q $ --- to the exotic solution of GYBE for 
$ \tilde{GL}_q(1|1)$ \cite{A-93}.

\section{ Quantum planes  and   groups } 

The standard solution of QYBE is given by (\ref{5}) with
$ \zeta=q. $ The generators of  $GL_q(2)$ have the commutation relations \cite{FRT-89}
\begin{equation}
 T=\left(\begin{array}{cc}
a & b \\
c &d \end{array}\right),
ab=qba, ac=qca, bd=qdb, cd=qdc, bc=cb, [a,d]=\lambda bc
\label{6}  
  \end{equation} 
and corresponding  quantum vector plane $C_q(2)$ is generated by $X^t=(x,y)^t$ such that $xy=qyx.$

In the standard for quantum groups basis (\ref{6}) in commutative limit $(q=1)$ the invariant $inv=X^tC_0X$ is
given by the matrix $(C_0)_{ik}=\delta_{i,3-i}$ with the units on the secondary diagonal, whereas in the really standard 
in many physical applications
Cartesian basis the invariant $inv=X^tIX$ is connected with 
the unit matrix $I.$ 
Therefore first of all we perform the linear transformation of basis 
of quantum plane $C_q(2)$ to the Cartesian one
\begin{equation}
X=DY, \quad
\left(\begin{array}{c}
x \\
y \end{array} \right)=
\frac{1}{\sqrt{2}}
\left(\begin{array}{cc}
1 & -i \\
1 &i \end{array}\right)
\left(\begin{array}{c}
p \\
r \end{array} \right),
 \label{7}  
  \end{equation} 
which induce the similarity transformation $U=D^{-1}TD $ of the generators of the quantum group $GL_q(2).$ Commutation relation of Cartesian generators of $C_q(2)$ are now in the form
\begin{equation} 
[r,p]=i\frac{q-1}{q+1}(r^2+p^2)=i(r^2+p^2)\tanh \frac{z}{2},  
\label{8}  
  \end{equation} 
 where $q=\exp z.$
  
  Let us substitute $ \iota \hat{r} $ instead of $r$ and $ \iota v $
 instead of $z,$ i.e. $\exp \iota v =1 + \iota v,$ 
 where $\iota $ is nilpotent unit $ \iota^2=0.$ Then the commutation relation for the new generators is as follows: 
 $[\hat{r},p]=hp^2, h=i\frac{v}{2}$ and is coincide with that  of 
 $h$-plane $C_h(2)$ \cite{AKS-94}. In the commutative case the above
 substitution $ \iota \hat{r} $ correspond to the transition from Euclidean
 plane to the flag (fiber) plane \cite{G-90}. Therefore, {\it $h$-plane $C_h(2)$ is the non-commutative (quantum) analog of the flag plane.}
 The associated quantum group $GL_h(2)$ is obtained from $GL_q(2)$
 by contraction,  their Cartesian generators are subject of the commutation relations
$$
U(\iota)=\left(\begin{array}{cc}
s & \iota t \\
\iota u & w \end{array}\right),\quad
 [s,w]=0,\quad [u,t]=h(s+w)(t+u),
 $$
 \begin{equation}
  [s,t]=[u,s]=hs(s-w),\quad
 [t,w]=[w,u]=hw(s-w) 
\label{9}  
  \end{equation}
and the following formulae   
  \begin{equation} 
 \delta\left(\begin{array}{c}
p \\
\iota \hat{r}\end{array}\right)=
\left(\begin{array}{cc}
s & \iota t\\
\iota u & w \end{array}\right)\dot{\otimes}
\left(\begin{array}{c}
p \\
\iota \hat{r}\end{array}\right)=\left(\begin{array}{c}
s\otimes p \\
\iota(u\otimes p+w\otimes \hat{r}) \end{array}\right) 
\label{10}  
  \end{equation} 
describe co-action of $GL_h(2)$ on $C_h(2).$

The exotic solution of QYBE is given by (\ref{5}) with
$ \zeta=-q^{-1}. $ The generators of the corresponding quantum plane
 $D_q(2)$ commute in the following way \cite{H-91}
\begin{equation} 
 xy=qyx, \quad  y^2=0. 
\label{11}  
  \end{equation} 
Dual (or Study) numbers $x+\iota\tilde{y}, \;  x,\tilde{y} \in R$
may be defined by the relations: $ xy=yx, \;   y^2=0,$ where 
$y=\iota\tilde{y}.$ Similarly to the complex plane the dual numbers
form the dual plane $D(2).$
 Therefore, {\it  quantum plane $D_q(2)$ is the non-commutative (quantum) analog of the dual plane $D(2).$} 
 If one introduce the generating matrix $T$ of the associated 
quantum group $\tilde{GL}_q(2)$ as in (\ref{6}), then their commutation
relations may be written as
\begin{equation} 
 b^2=c^2=0,\; bc=cb,\; ac=qca,\; db=-qbd,\; dc=-qcd,\; [a,d]=\lambda bc.  
\label{12}  
  \end{equation} 
 Let us stress that all generators of $D_q(2)$ and $\tilde{GL}_q(2)$
 are of the even order and relations $y^2=b^2=c^2=0 $ are appear due to
 their nilpotent nature.

\section{ Quantum superplanes  and  supergroups} 

The standard solution of GYBE is given by (\ref{5}) with
$ \zeta=q^{-1}. $ The generators of quantum supergroup $GL_q(1|1)$ have  commutation relations 
$$  
  T=\left(\begin{array}{cc}
a & \alpha \\
\beta &b \end{array}\right), \quad
\alpha^2=\beta^2=0,\;\; \{\alpha,\beta \}=0, 
$$
\begin{equation}
a\alpha=q\alpha a,\;\; a\beta=q\beta a,\;\;
 b\alpha=q\alpha b,\;\;  b\beta=q\beta b,\;\; [a,b]=\lambda \beta \alpha
\label{13}  
  \end{equation} 
and thouse for the corresponding quantum superplane $C_q(1|1)$ are
$x\theta=q\theta x, \; \theta^2=0. $ Generators $a,b,x$ are even  and
$\alpha, \beta, \theta $ are odd (Grassmann). Performing the superlinear transformation
of generators
\begin{equation} 
 x=y+\frac{h}{v}\xi, \quad \theta =\xi +\frac{h}{v}y, 
\label{14}  
  \end{equation} 
where $h$ is odd, $ h^2=0$ and $q=1+v, $ we obtain  $C_q(1|1)$ in the new  basis $y, \xi:$
\begin{equation} 
 [y,\xi]=hy^2+v\xi y, \quad \xi^2=-h\xi y. 
\label{15}  
  \end{equation} 
  After contraction $v\rightarrow 0$ the quantum $h$-superplane
$C_h(1|1)$ \cite{DP-95} is achieved
\begin{equation} 
 C_h(1|1)=\{y, \xi| [y,\xi]=hy^2, \quad \xi^2=-h\xi y \}, 
\label{16}  
  \end{equation} 
 which may be interpreted as {\it non-commutative (quantum) analog
 of the flag superplane} \cite{GKK-03}. The quantum supergroup
 $GL_h(1|1):$ 
  $$  
  T=\left(\begin{array}{cc}
m & \psi \\
\varphi &n \end{array}\right), \quad
\psi^2=0,\;\; \varphi^2=h\varphi(n-m), \;\; \{\psi ,\varphi\}=h\psi (n-m), 
  $$
  $$
 [m,\varphi ]=h(\varphi\psi -m(n-m)),\;\;
 [n,\varphi ]=h(\varphi\psi -n(n-m)),  
  $$
\begin{equation}
[m,\psi ]=[n,\psi ]=0,\;\; [n,m]=h\psi(n-m)
\label{17}  
  \end{equation} 
 co-act on $C_h(1|1)$ according with (\ref{4}). 
  
 The exotic solution of GYBE is given by (\ref{5}) with
$ \zeta=-q. $ The  corresponding quantum graded plane is
$\tilde{C}_q(1|1)=\{z,\mu |\; z\mu =q\mu z \} $
(but no relation $\mu^2=0$) and the quantum graded group
$\tilde{GL}_q(1|1)$ is given by 
$$  
  T=\left(\begin{array}{cc}
c & \gamma \\
\delta &d \end{array}\right), \quad
\{\delta,\gamma \}=0,\;\; c\gamma =q\gamma c, \;\; c\delta =q\delta c, 
  $$
  \begin{equation}
\gamma d=-qd\gamma,\;\;
 \delta d =-qd\delta, \;\;
[c,d]=\lambda \delta \gamma,
\label{18}  
  \end{equation}
but no relations  $ \gamma^2=\delta^2=0. $
After the similar to (\ref{14})  transformation
\begin{equation} 
 z=t+\frac{h}{v}\nu, \quad \mu =\nu +\frac{h}{v}t, 
\label{19}  
  \end{equation} 
where $h$ is odd, $ h^2=0, q=1+v, $ generator $t$ is even,
$\nu$ is odd and by hand $\nu^2=0,$  we obtain
$\tilde{C}_q(1|1)=\{t,\nu |\; [t,\nu ] =ht^2 + v\nu t, \;\; \nu^2=0 \}. $
Contraction $v\rightarrow 0 $ gives in result {\it new exotic $h$-superplane}
$\tilde{C}_h(1|1)=\{t,\nu |\;[t,\nu ] =ht^2, \; \nu^2=0 \}.$

\section{Quantum superspace $C_q(1|2)$ and its contraction} 

The standard  $N=3$ solution of GYBE associated with the quantum supergroup
$GL_q(1|2)$ and quantum superspace $C_q(1|2)$ is in the form \cite{I-95}
\begin{equation}
R_q=\left(\begin{array}{ccccccccc}
q& 0 & 0       & 0 & 0 & 0        & 0 & 0 & 0\\
0& 1 & 0       & 0 & 0 & 0        & 0 & 0 & 0\\
0& 0 & 1       & 0 & 0 & 0        & 0 & 0 & 0\\

0&\lambda & 0  & 1 & 0 & 0        & 0 & 0 & 0 \\ 
0& 0 & 0       & 0 & q^{-1} &0    & 0 & 0 & 0 \\
0& 0 & 0       & 0 & 0 & 1        & 0 & 0 & 0 \\

0& 0 & \lambda & 0 & 0 & 0        & 1 & 0 & 0 \\
0& 0 & 0       & 0 & 0 &-\lambda  & 0 & 1 & 0 \\
0& 0 & 0       & 0 & 0 & 0        & 0 & 0 &q^{-1} 
\end{array} \right).  
   \label{20}  
  \end{equation} 
According with (\ref{3}) the  quantum superspace $C_q(1|2)$ is generated by even $x$
and odd  $\theta_1, \theta_2$ with commutation relations \cite{I-95}
\begin{equation} 
 x\theta_k=q\theta_kx, \;\;
\theta_1\theta_2 =-q\theta_2\theta_1, \;\;
\theta_k^2=0, \;\; k=1,2. 
\label{21}  
  \end{equation} 
After the linear transformation of the  generators
 \begin{equation}
 x=x, \quad 
 \theta_1 ={\frac{1}{\sqrt{2}}} (\xi_1 -i\xi_2), \quad
 \theta_2 ={\frac{1}{\sqrt{2}}} (\xi_1 +i\xi_2)
\label{22}  
  \end{equation}
quantum superspace in  Cartesian basis looks as follows
\begin{equation} 
 C_q(1|2)=\{x,\xi_k |\;  x\xi_k=q\xi_k,\; k=1,2,\; \{\xi_1,\xi_2\}=0, \; \xi_1^2=\xi_2^2=i{\frac{q-1}{q+1}} \xi_1\xi_2 \}. 
 \label{23}  
  \end{equation} 
  
 A superspace with flag (fiber) odd subspace is obtained \cite{GKK-03}  by substitution $\xi_2 \rightarrow \iota\xi_2 $ and  transformation $q=e^z \rightarrow e^{\iota v}$ of deformation parameter need be added
 in quantum case. As the result we have  quantum flag superspace 
 \begin{equation} 
  C_h(1|2;\iota)=\{x,\xi_k |\; [x,\xi_k]=0, \;  \{\xi_1,\xi_2\}=0,\;
  \xi_1^2=0, \; \xi_2^2=h\xi_1\xi_2, \;h=i\frac{v}{2} \}
 \label{24}  
  \end{equation} 
 with the fiber odd subspace.  Generator $x$ commute with $\xi_1, \xi_2,$
 therefore we may put $x=1$ and the quotient algebra
  \begin{equation}
 \hat{C}_h(2)=C_h(1|2;\iota)/\{x=1\}=
 \{\xi_1,\xi_2 |\;
 \{\xi_1,\xi_2\}=0,\;\xi_1^2=0, \;\xi_2^2=h\xi_1\xi_2 \}
 \label{25}  
  \end{equation}
 is geometrically interpreted as {\it the flag superplane with both odd generators.} It is interesting to note that commutation relations (\ref{25}) are the same as for the differentials
 of $h$-plane $C_h(2)$ (see \cite{ACDM-01}, (6.11--14)).
  Co-action of consistent with $ C_h(1|2;\iota)$  quantum group 
  $SL_h(1|2;\iota )$ is given by
 \begin{equation}
U(\iota)\dot{\otimes}X=\left(\begin{array}{ccc}
1& 0 & 0 \\
0 & k & \iota r \\
0 & \iota m & k^{-1} \end{array} \right)\dot{\otimes}
\left(\begin{array}{c}
1 \\
\xi_1 \\
\iota\xi_2 \end{array}\right)=
\left(\begin{array}{l}
1\otimes 1 \\
k\otimes \xi_1 \\
\iota(m\otimes\xi_1+k^{-1}\otimes\xi_2)\end{array}\right),
\label{26}  
  \end{equation}
where generators of $SL_h(1|2;\iota )$ are subject of commutation relations
$$ 
[r,k]=[k,m]=h(k^2-1),\quad
[k^{-1},r]=[m,k^{-1}]=h(1-k^{-2}),
$$
\begin{equation}
[r,m]=h(k+k^{-1})(r+m).
\label{27}  
  \end{equation}

\section{Conclusion}

We have investigated quantum (super)groups and quantum (super)planes,
which correspond to the  simplest one-parameter $GL$-type solutions of YBE from
geometrical point of view. After linear transformation of generators to the Cartesian
basis and standard contraction the well-known $h$-plane 
$C_h(2)=\{ [\hat{r},p]=hp^2 \}$ with both even generators is interpreted
as the non-commutative flag plane. Quantum plane (\ref{11})
$ D_q(2)= \{ xy=qyx, \;  y^2 \},$ associated with exotic solution of QYBE,
may be regarded as the non-commutative   dual plane.
  The standard and exotic solutions of GYBE lead after contractions to
the non-commutative flag superplane (\ref{16}) 
$C_h(1|1)=\{[y,\xi]=hy^2, \; \xi^2=-h\xi y \}$
and new exotic flag superplane 
$\tilde{C}_q(1|1)=\{ [t,\nu ] =ht^2, \; \nu^2=0 \} $ 
with one even and one odd generators.
New non-commutative flag superplane (\ref{25})
$ \hat{C}_h(2)= \{ \{\xi_1,\xi_2\}=0,\;\xi_1^2=0, \;\xi_2^2=h\xi_1\xi_2 \}$
  with both odd generators is obtained from $C_q(1|2)$ by the linear transformation of the odd generators, contraction and quotient on even  generator.

\end {document}